\documentclass[reqno, 12pt]{amsart}
\usepackage{amsfonts}
\usepackage{amssymb,amsmath}
\usepackage{amsthm}
\usepackage{upref}
\usepackage{enumerate}
\usepackage{bbm}

\usepackage{amscd}

\makeatletter
\def\LaTeX{\leavevmode L\raise.42ex
    \hbox{\kern-.3em\size{\sf@size}{0pt}\selectfont A}\kern-.15em\TeX}
 \makeatother

\numberwithin{equation}{section}

\newtheorem{lemma}{Lemma}[section]
\newtheorem{theorem}[lemma]{Theorem}

\theoremstyle{definition}

\newtheorem{example}[lemma]{Example}

\newtheorem{remark}[lemma]{Remark}

 \newcommand{\supp}{\operatorname{supp}}
  \newcommand{\e}{\eqref}

\newcommand{\q}{\quad}

\newcommand{\wt}{\widetilde}

\newcommand{\z}{\zeta}
\newcommand{\ov}{\overline}
 \renewcommand{\d}{\delta}
 
 \newcommand{\cl}{\operatorname{clos}}

\renewcommand\Im{\operatorname{Im}}

\def\qqq{\mathrel{\subset\mkern-15mu\lower.38ex\hbox{${\scriptscriptstyle\rightarrow}$}}}

\let\cal\mathcal

\let\Bbb\mathbb

\begin{document}

\title [   Toeplitz {\tiny versus}  Hankel] {      Toeplitz versus  Hankel: semibounded operators  }
 
\author{ D. R. Yafaev}
\address{ IRMAR, Universit\'{e} de Rennes I, Campus de
  Beaulieu,  Rennes, 35042  FRANCE and SPGU, Univ. Nab. 7/9, Saint Petersburg, 199034 RUSSIA}
\email{yafaev@univ-rennes1.fr}


\thanks {Supported by  project   Russian Science Foundation   17-11-01126}

\begin{abstract}
Our goal is to compare various results for
Toeplitz $T$ and Hankel $H$ operators.  We consider semibounded operators and find necessary and sufficient conditions for their quadratic forms to be closable. This property allows one to define 
$T$ and   $H$ as self-adjoint operators under minimal assumptions on their matrix elements. We also describe domains of the closed Toeplitz and Hankel quadratic forms.
   \end{abstract}

\maketitle


\section{Introduction. Bounded operators}  

{\bf 1.1.}
This  is a short survey based on the talk given by  the  author at the conference "Spectral Theory and Applications"   held in May 2017 in Krakow. Our aim is to compare various properties of Hankel  and Toeplitz   operators.
 We refer to the books \cite{Bo, GGK, NK, Pe} for basic information on  
these classes of operators.

Formally,    Toeplitz $T$ and Hankel  $H$ operators      are defined in the space  $\ell^2  ({\Bbb Z_{+}})$  of sequences $f=(f_{0}, f_{1}, \ldots)$ by the relations 
   \begin{equation}
(T f )_{n}=\sum_{m\in {\Bbb Z}_{+}} t_{n-m} f_{m} , \q n\in {\Bbb Z}_{+},
\label{eq:Toe}\end{equation}
and 
  \begin{equation}
(H f )_{n}=\sum_{m\in {\Bbb Z}_{+}} h_{n+m} f_{m} , \q n\in {\Bbb Z}_{+}.
\label{eq:B1h}\end{equation}

Let us also introduce   discrete convolutions in the space  $\ell^2  ({\Bbb Z}) $ (known also as  Laurent operators) acting   by the formula
\begin{equation}
(L g )_{n}=\sum_{m\in {\Bbb Z} } t_{n-m} g_{m}, \q  g= \{ g_n\}_{n\in{\Bbb Z} }.
\label{eq:Lor}\end{equation}
Of course, by the discrete Fourier transform, $L$ reduces to the operator of multiplication by the function (symbol)
\begin{equation}
t(z)=\sum_{n\in {\Bbb Z} } t_{n} z^n,
\label{eq:st}\end{equation}
and so its spectral analysis is trivial.
If the sequences  $t_{n}$ in \e{eq:Toe} and \e{eq:Lor} are the same, one might  expect that properties of  the  operators $T$ and $L$ are also similar.  Surprisingly, this very naive conjecture is not totally wrong.

 The precise definitions of the operators   $T$  and $H$ require some accuracy.  Let  $\cal D \subset \ell^2({\Bbb Z}_{+})  $ be the dense set   of sequences $f=\{f_{n}\}_{n\in {\Bbb Z}_{+}} $ with only a finite number of non-zero components. If the sequences 
 $t=\{t_{n}\}_{n\in {\Bbb Z}_{+}}  \in \ell^2({\Bbb Z}_{+})  $ and  $h=\{h_{n}\}_{n\in {\Bbb Z}_{+}}  \in \ell^2({\Bbb Z}_{+})  $, then for $f\in \cal D $, the vectors $Tf$ and $Hf$
    belong to $ \ell^2({\Bbb Z}_{+})  $ so that the operators $T$ and $H$ are defined on $\cal D $.
 Without such a priori assumptions, 
  only   Toeplitz  
   \begin{equation}
t[f,f] =\sum_{n,m\geq 0} t_{n-m} f_{m}\bar{f}_{n}   
 \label{eq:QFq}\end{equation}
 and
 Hankel  
        \begin{equation}
h[f,f] =\sum_{n,m\geq 0} h_{n+m} f_{m}\bar{f}_{n}   
 \label{eq:Ha}\end{equation}
  quadratic forms
are well defined for all $f\in \cal D$.

  Obviously, a Toeplitz operator $T$ (resp., a Hankel operator $ H$) is bounded if and only if the estimate
   \begin{equation}
|t[f,f]| \leq C  \| f\|^2 ,\q f\in \cal D,
 \label{eq:QFT}\end{equation}
 or
    \begin{equation}
|h[f,f]| \leq C  \| f\|^2, \q f\in \cal D,
 \label{eq:QFH}\end{equation}
 is satisfied. Here and below $\| f\|$ is the norm of $f$ in the space $\ell^2({\Bbb Z}_{+}) $; $C$ are different positive constants; $I$ is the identity operator.
 
 \medskip

{\bf 1.2.}
   Let us recall a necessary and sufficient condition  for   Toeplitz and Hankel operators   to be bounded.
   In terms of quadratic forms the conditions of boundedness of these operators   can be stated without any a  priori assumptions on their matrix elements.
    Below $d{\bf m}(z)=(2\pi i z)^{-1} dz$ is   the normalized Lebesgue measure on  the unit circle $\Bbb T$. For $p\geq 1$, we set 
   $L^p (\Bbb T ) =L^p (\Bbb T ; d{\bf m})$.

 \begin{theorem}[Toeplitz]\label{B-H} 
 Estimate \e{eq:QFT} is true if and only if the $t_{n}$ are the Fourier coefficients of some bounded function $t (z) $ on  $\Bbb T$:
    \begin{equation}
  t_{n} =\int_{\Bbb T} z^{-n} t (z) d{\bf m} (z) , \q n\in{\Bbb Z}, \q t \in L^\infty (\Bbb T ).
 \label{eq:BH}\end{equation}
     \end{theorem}
     
  Thus a Toeplitz operator $T$ is bounded if and only if the corresponding Laurent operator \e{eq:Lor} is bounded.
 
 The following result is due to Z.~Nehari \cite{Nehari}.

     \begin{theorem}\label{Neh} 
 Estimate \e{eq:QFH} is true if and only if there exists
  a bounded function $h(z)$ on  $\Bbb T$ such that  
    \begin{equation}
  h_{n} =\int_{\Bbb T} z^{-n} h (z) d{\bf m} (z) , \q n\in{\Bbb Z}_{+}, \q h \in L^\infty (\Bbb T ) .
 \label{eq:BHN}\end{equation} 
     \end{theorem}
     
     Despite a formal similarity, Theorems~\ref{B-H} and \ref{Neh} are essentially different because the symbol $t(z)$ of a Toeplitz operator $T$ is uniquely defined by 
     relation  \e{eq:BH} while   \e{eq:BHN} imposes conditions only on the Fourier coefficients $h_{n}$ of $h(z)$ with $n\in{\Bbb Z}_{+}$. So among the functions satisfying \e{eq:BHN} there may be both bounded and unbounded functions. The following example illustrates this phenomenon.
     
     \begin{example}\label{Hilb} 
     Let $ h_{n}= (n+1)^{-1}$
     (the corresponding Hankel operator $H$ is known as the Hilbert matrix). The ``natural" symbol
     \[
     h(z)=\sum_{n\in{\Bbb Z}_{+}} (n+1)^{-1} z^n
     \]
     is unbounded on $\Bbb T$ (at the single point $z=1$). However the function
     \[
     \wt{h}(z)=1+\sum_{n\geq 1} (n+1)^{-1} (z^n -z^{-n})
     \]
     is also a symbol of $H$ and  $ \wt{h}\in L^\infty (\Bbb T)$. Therefore $H$ is a bounded operator.
     \end{example}
      

Actually,  properties of Toeplitz and Hankel operators are quite different. For example,   a Toeplitz operator $T$ is never compact unless $T=0$ (see, e.g., Section~3.1 in the book \cite{Pe}. On the contrary, a Hankel operator $H$ is   compact if its symbol can be chosen as a continuous function. Properties of compact Hankel operators are very thoroughly studied in \cite{Pe}.

 \medskip

{\bf 1.3.}
The  results about unbounded   operators     are very scarce.  We can mention only the paper \cite{Hartman} by P.~Hartman and the relatively recent survey \cite{Sarason} by D.~Sarason; see also references in these articles. These articles are devoted to Toeplitz   operators.
   Note that the theory of general unbounded integral operators was initiated by T.~Carleman in \cite{Carleman}, but this theory   does not practically provide concrete conditions guaranteeing, for example, that  a given symmetric operator is essentially self-adjoint. Probably, it is impossible to develop a complete theory for general integral operators, but it is very tempting to do this for Toeplitz and Hankel operators   possessing special structures.

  Our goal here is to describe exhaustive results in the semibounded case, both for Toeplitz and Hankel operators. It looks instructive to compare the results for these two very different classes.
    Our approach relies on a certain
  auxiliary algebraic construction combined with some classical analytical results. The  algebraic construction  is more or less the same for Toeplitz and Hankel quadratic forms, but the analytic results we use are quite different. In particular, the Riesz Brothers theorem plays the crucial role for Toeplitz   operators,  while  the Paley-Wiener  theorem is important for 
  Hankel operators.

In Section~2, we find necessary and sufficient conditions for Toeplitz and Hankel quadratic forms to be closable.
In Section~3, we describe their closures. Finally, in Section~4, we very briefly (see  \cite{YaWH}, for details) discuss Wiener-Hopf operators that are a continuous analogue of Toeplitz   operators. To treat Wiener-Hopf operators,  we need a continuous version of the Riesz Brothers theorem which was not available in the literature.
  
\section{Semibounded operators and their quadratic forms}  

{\bf 2.1.}
In the semibounded case, it is natural to define operators via their quadratic forms. The corresponding construction has abstract nature. It  is due to Friedrichs and is described, for example, in the book \cite{BSbook}.  Consider 
an arbitrary Hilbert space $\cal H$  with the norm $\| \cdot \|$ and a
real quadratic form $b[f,f]$   defined on a set $\cal D$ dense in   $\cal H$. Assume that the form $b[f,f]$ is semibounded, that is, 
\begin{equation}
b[f,f] \geq \gamma \| f\|^2  ,\q f \in{\cal D},  
 \label{eq:T1}\end{equation}
 for some  $\gamma\in \Bbb R$. Suppose  first that $\gamma>0$. Then one can introduce a new norm
 \[
 \|f\|_{b} =\sqrt{b[f,f]}
 \]
 which is stronger than the initial norm $\|f\|$. If $\cal D$ is a complete Hilbert space  with respect to the norm
 $\| \cdot \|_{b}$ (in this case the form $b$ is called closed), then there exists a unique self-adjoint operator $B$ in $\cal H$ with the domain ${\cal D} (B)\subset {\cal D}$ such that $B\geq \gamma I$ and
        \[
        b[f, g] =  ( Bf, g), \q \forall f \in {\cal D}(B), \q \forall g\in  {\cal D} .
       \]
      Moreover,  ${\cal D}(\sqrt{B}) = {\cal D}$ and 
      \[
        b[f,f]= \| \sqrt{B} f\|^2, \q \forall f\in  {\cal D} .
        \]
Thus the self-adjoint operator $B$ is correctly defined although       its domain $ {\cal D} (B) $   does not admit an efficient description.  
              
              If $\cal D$ is not $b$-complete, then of course one can take its  completion $\cal D [b]$ in the norm $\| \cdot \|_{b}$, extend $  b[f,f]$ by continuity onto $\cal D [b]$ and then try to apply the construction above to the   form  $  b[f,f]$ defined on $\cal D [b]$. However this procedure meets with an obstruction because, in general, $\cal D [b]$ cannot be realized as a subset of $\cal H$. One can avoid this problem only for the so-called closable forms. By definition, a form $b[f,f]$ is closable if the conditions 
   \[
\| f^{(k)}\|\to 0 \q{\rm and}  \q \| f^{(k)} - f^{(j)}\|_b\to 0
\]
 as $k,j\to\infty$ imply that $\| f^{(k)}\|_{b}\to 0$. In this case $\cal D [b]\subset\cal H$, the form $  b[f,f]$ is closed on $\cal D [b]$,  and so there exists a self-adjoint operator $B$ corresponding to this form.

 If $\gamma$ in    \e{eq:T1} is not positive, then one applies the definitions above to a form $  b_\beta [f,f]= b[f,f] + \beta \| f\|^2$ for some $\beta>  -\gamma$ and   defines $B$ by the equality $B=B_\beta - \beta I$. So,  we can suppose that the number $\gamma$ in \e{eq:T1} is positive; for definiteness, we choose $\gamma=1$. 
 
 To summarize,    semibounded self-adjoint operators are correctly defined if and only if the corresponding quadratic forms are closable.
  Of course not all forms are closable. On the other hand, it is easy to see that if $B_{0}$ is a symmetric semibounded operator with domain ${\cal D} (B_{0})$, then the form $  b[f,f]=  (  B_{0} f,f)$ defined on   ${\cal D} (B_{0})$ is necessarily closable.

  \medskip

{\bf 2.2.}
Let us come back to Toeplitz and Hankel operators when   ${\cal H}  = \ell^2 ({\Bbb Z}_{+})$ and 
 $\cal D \subset \ell^2({\Bbb Z}_{+})  $ consists of sequences $f=\{f_{n}\}_{n\in {\Bbb Z}_{+}} $ with only a finite number of non-zero components. We now suppose that $t_{n}= \bar{t}_{-n}$ and $h_{n}= \bar{h}_{n}$ for all $n\in{\Bbb Z}_{+}$.    
  First, we state the conditions for Toeplitz and Hankel quadratic forms to be semibounded.
 
 For Toeplitz   quadratic forms, we use  the following well known result (see, e.g., \S 5.1 of the book \cite{AKH}) that is a consequence of the F.~Riesz-Herglotz theorem.

  \begin{theorem}\label{R-H} 
   The condition 
      \[
 \sum_{n,m\geq 0} t_{n-m} f_{m}\bar{f}_{n}\geq 0, \q \forall f\in \cal D,
 \]
  is satisfied if and only if there exists a
  non-negative $($finite$)$ measure $dM(z)$ on the unit circle $\Bbb T$ such that the coefficients $t_{n }$ admit the representations 
  \begin{equation}
  t_{n} =\int_{\Bbb T} z^{-n} dM (z)  , \q n\in{\Bbb Z}.
 \label{eq:WH}\end{equation}
    \end{theorem}
    
    Equations \e{eq:WH} for the measure $dM (z)$ are known as the trigonometric moment problem. 
       Of course  their solution is   unique.  Note that the identity   $I$ is the Toeplitz operator (with $t_{0 }=1$ and $t_{n }=0$ for $n\neq 0$) and the corresponding measure $dM (z)$ in \e{eq:WH} is  the normalized Lebesgue measure $d {\bf m} (z) $. Therefore the measure corresponding to the form
    $t[g,g] + \beta \| g\|^2 $ equals $dM  (z)+ \beta d {\bf m} (z)$. So we have a one-to-one correspondence between Toeplitz  quadratic forms satisfying estimate  \e{eq:T1} and real measures satisfying the condition $M  (X)\geq \gamma    {\bf m} (X)$ for all Borelian sets $X\subset{\Bbb T}$.
     
      Hankel quadratic forms are linked to the  power moment problem. The following result obtained by Hamburger  in \cite{Hamb} plays the role of Theorem~\ref{R-H}.

   \begin{theorem}\label{Hamb} 
   The condition 
      \begin{equation}
 \sum_{n,m\geq 0} h_{n+m} f_{m}\bar{f}_{n}\geq 0, \q \forall f\in \cal D,
 \label{eq:QFh}\end{equation}
  is satisfied if and only if there exists a
  non-negative measure $d{\sf M}(x)$ on $\Bbb R$ such that the coefficients $h_{n }$ admit the representations 
    \begin{equation}
h_{n } = \int_{-\infty}^\infty x^n d {\sf M}(x), \q \forall n=0,1,\ldots.
   \label{eq:WHh}\end{equation}
   \end{theorem}
    
    \medskip
    

{\bf 2.3.}
If conditions 
 \begin{equation}
 t[f,f] \geq \gamma \| f\|^2, \q \forall f\in \cal D,
 \label{eq:QTT}\end{equation}
 (for some $\gamma\in {\Bbb R}$) or \e{eq:QFh} are satisfied and 
     \begin{equation}
\sum_{n\in\Bbb Z}|t_{n}|^2 <\infty \q {\rm or}  \q \sum_{n\in {\Bbb Z}_{+}} h_{n}^2 <\infty ,
\label{eq:Tx2}\end{equation}
 then the   forms  $t[f,f]$ or $h[f,f]$  are   closable. Indeed, in this case  the Toeplitz operator \e{eq:Toe} or Hankel operator \e{eq:B1h} are well defined and symmetric
 on the set ${\cal D}  $.
However, the conditions \e{eq:Tx2} are by no means necessary  for the forms $t[f,f]$ or  $h[f,f]$  to be  closable.

   Our main goal is to find {\it necessary and sufficient} conditions for the forms $t[f,f]$ and $h[f,f]$  to be  closable. The answers to these questions are strikingly  simple.
   
    We start with Toeplitz forms.
   
    \begin{theorem}\label{T1}\cite[Theorem 1.3]{YaT}  
    Let the form $t[f,f]$ be given by formula \e{eq:QFq} on elements $f\in {\cal D}$, and let the condition 
 \e{eq:QTT} be satisfied. 
   Then the form $t[f,f] $ is closable in the space $\ell^2 ({\Bbb Z}_{+})$ if and only if the measure $dM  (z)$ in   the equations \e{eq:WH} is absolutely continuous.
    \end{theorem}
    
     \begin{example}\label{Tex} 
      If $t_{n}=1$ for all $n\in {\Bbb Z}$, then $M (\{1\})=1$ and $M ({\Bbb T} \setminus \{1\})=0$. This measure
    is supported by the single point $z=1  $, and the corresponding quadratic form is not closable.  
    \end{example}
    
    Of course  Theorem~\ref{T1} means that $dM  (z)=  t(z) d{\bf m}(z)$ where the function $t \in L^1 (\Bbb T; d{\bf m})$ and $t (z)\geq \gamma$. Thus Theorem~\ref{T1} extends Theorem~\ref{B-H} to semibounded operators. The function 
$t (z)$ is known as the symbol of the Toeplitz operator     $T$. So, Theorem~\ref{T1}  shows that for a semibounded Toeplitz operator (even defined via the corresponding quadratic form), the symbol exists and is a semibounded function.

The result for Hankel quadratic forms is stated as follows.

   \begin{theorem}\label{Hamb1}\cite[Theorem 1.2]{Yunb}   
     Let   assumption \e{eq:QFh} be satisfied.  Then the following conditions are equivalent: 
     
     \begin{enumerate}[\rm(i)]
\item
The form $h[f,f]$ defined on $\cal D$  is closable   in the space $\ell^2 ({\Bbb Z}_{+})$.
 
\item
The measure $d {\sf M}  (x)$ defined by equations  \e{eq:WHh} satisfies  the condition
  \begin{equation}
{\sf M} ({\Bbb R}\setminus (-1,1) )=0
\label{eq:dense}\end{equation} 
$($to put it differently, $\supp {\sf M}\subset [-1,1]$ and ${\sf M}(\{-1\}) = {\sf M}(\{1\})=0)$.

 \item
The matrix elements  $ h_n\to 0$ as $n\to \infty$.
\end{enumerate}
        \end{theorem}
        
       \begin{remark}\label{Hex} 
  In general, the measure $d{\sf M}  (x)$ in \e{eq:WHh} is not unique, but it is unique under the assumptions of Theorem~\ref{Hamb1}.
    \end{remark}  
 
          Theorem~\ref{Hamb1}  is to a large extent motivated by the following classical result of H.~Widom.
           
              \begin{theorem}\label{Widom}\cite[Theorem~3.1]{Widom}
                Let   the matrix elements  $ h_n $ of the Hankel operator \e{eq:B1h}  be given by the equations 
                   \[
h_{n } = \int_{-1}^1 x^n d {\sf M} (x), \q \forall n=0,1,\ldots,  
 \]
  with some non-negative  measure $dM  (x)$.
            Then the following conditions are equivalent: 
                
\begin{enumerate}[\rm(i)]
\item
The operator $H$ is bounded.   
   \item
 ${\sf M} ((1-\varepsilon,1])= O(\varepsilon)$ and  ${\sf M}([-1, -1+\varepsilon))= O(\varepsilon)$ as $\varepsilon\to 0$.
  \item
  $h_{n}=O(n^{-1})$ as $n\to \infty$.
\end{enumerate}
    \end{theorem}

 We emphasize that, in the semibounded case,   Theorems~\ref{T1} and \ref{Hamb1}  give optimal conditions for Toeplitz and Hankel operators to be defined as self-adjoint operators. Below we briefly discuss the proofs of these results.

 \medskip




  \medskip
  
  {\bf 2.4.}  
We start with  Theorem~\ref{T1}, where we may suppose that $\gamma=1$ in  \e{eq:QTT}. Set  
   \begin{equation}
( {\bf A} f) (z)=\sum_{n=0}^\infty   f_{n} z^n.
 \label{eq:A}\end{equation}
 Then  
  \begin{equation}
   \| {\bf A} f \|_{ L^2 (\Bbb T )} =  \| f\|_{\ell^2 ({\Bbb Z}_{+})}
 \label{eq:A2P}\end{equation}
 for all $f\in \ell^2 ({\Bbb Z}_{+})$.
 Clearly, ${\bf A}$ is a unitary mapping of  $  \ell^2 ({\Bbb Z}_{+})$  onto the Hardy space   $H^2 ({\Bbb T} )$ of functions analytic in the unit disc.   In view of  equations  \e{eq:WH}, we also have
   \begin{equation}
      \| {\bf A} f \|^2_{ L^2 (\Bbb T; dM)} =t [f,f] , \q f\in {\cal D}.
 \label{eq:A2}\end{equation}

  The `` if " part of  Theorem~\ref{T1} is quite easy. Suppose  that for a sequence $f^{(k)}\in \cal D$
 \[
 \| f^{(k)}\|_{\ell^2 ({\Bbb Z}_{+})}\to 0 \q {\rm and} \q  t [f^{(k)}-f^{(j)},f^{(k)}-f^{(j)}]\to 0
 \]
 as $k,j\to\infty$. Put $g^{(k)}={\bf A} f^{(k)}$. It follows from \e{eq:A2P}, \e{eq:A2}  that
   \begin{equation}
 \| g^{(k)}\|_{ L^2 (\Bbb T)}\to 0 \q {\rm and} \q    \| g^{(k)}-g^{(j)} \|^2_{ L^2 (\Bbb T; dM)} \to 0
 \label{eq:A2f}\end{equation}
 as $k,j\to\infty$. Since the space $L^2 (\Bbb T; dM)$ is  complete, there 
 exists a function $g\in L^2 (\Bbb T; dM)$ such that $g^{(k)}\to g$ in $L^2 (\Bbb T; dM)$ and hence in $L^2 (\Bbb T)$ . The first condition \e{eq:A2f} implies that
 $\| g \|_{ L^2 (\Bbb T)}=0$ whence  $\| g \|_{ L^2 (\Bbb T; dM)}=0$ because the measure $dM$ is absolutely continuous with respect to the Lebesgue measure. It now follows from \e{eq:A2} that 
 $$ t [f^{(k)},f^{(k)}]=\| g^{(k)} \|_{ L^2 (\Bbb T; dM)} \to 0$$
  as $k\to\infty$. Thus the form $t[f,f]$ is closable.

The proof of  the ``only if " part of  Theorem~\ref{T1} is less straightforward.
 Let us define  the operator $A \colon \ell^2 ({\Bbb Z}_{+})\to    L^2 (\Bbb T; dM)$ as the restriction of
  the operator $\bf A$ on the set $\cal D=:   \cal D (A)$.
 First, we note an assertion which is a  direct consequence of the identity \e{eq:A2}.

    \begin{lemma}\label{de}
 The form $   t [f,f]$ defined   on $\cal D$   is closable in the space $\ell^2 ({\Bbb Z}_{+})$  if and only if the operator $A  $  
  is closable.  
    \end{lemma}
     
      The next step is to construct the adjoint operator   $A^* \colon L^2 (\Bbb T; dM) \to \ell^2 ({\Bbb Z}_{+}) $.  Observe that       
 for an arbitrary $u\in L^2 (\Bbb T; dM)$, all the integrals
       \[
 \int_{\Bbb T}  u(z) z^{-n} dM(z)=: u_{n}, \q n\in {\Bbb Z}_{+},
\]
 are absolutely convergent and the sequence $ \{u_{n}\}_{n=0}^\infty$ is bounded. Let us
   distinguish a subset     ${\cal D}_{*}\subset   L^2 (\Bbb T; dM)$ by the condition   $\{u_{n}\}_{n=0}^\infty \in  \ell^2 ({\Bbb Z}_{+})$ for $u\in {\cal D}_{*}$.  For the proof of the following assertion, see  Lemma~2.4 in \cite{YaT}.

   \begin{lemma}\label{LTM}
 The operator $A ^*$ is given by the equality 
      \[
        (A^ {*}u)_{n}=  \int_{\Bbb T}  u(z) z^{-n} dM(z), \q n\in {\Bbb Z}_{+},
 \]
  on the domain ${\cal D}(A^*) ={\cal D}_{*}$. 
    \end{lemma}

      Recall that    an operator $A $   is closable if and only if its adjoint operator  $A^*$ is densely defined. We use the notation $ \cl{\cal D}_{*}$ for the closure of the set ${\cal D}_{*}$
      in the space $L^2 ({\Bbb T}; dM)$. Lemmas~\ref{de} and \ref{LTM}
      yield the following intermediary result.
      
       \begin{lemma}\label{adj}
 The operator $A  $ and the form $t[f,f]$ are closable if and only if  
   \begin{equation}
 \cl{\cal D}_{*} =L^2 ({\Bbb T}; dM). 
 \label{eq:D}\end{equation}
    \end{lemma}

Recall  the   Riesz Brothers theorem (see, e.g., Chapter~4 in \cite{Hof}) that we combine with the Parseval identity.

         \begin{theorem}\label{brothers}
       For a  complex $($finite$)$ measure $d\mu (z)$    on the unit circle $\Bbb T$, put
              \[
\mu_{n}=\int_{\Bbb T}  z^{-n}d\mu(z), \q n\in  {\Bbb Z},
\]
 and suppose that 
 \[
\sum_{n=0}^\infty |\mu_{n}|^2 <\infty. 
\]
     Then  the measure $d\mu   (z)$ is absolutely continuous.
    \end{theorem}

    We  need also   the following technical assertion (Lemma~2.7 in \cite{YaT}).
  
   \begin{lemma}\label{ac}
      Suppose that a set ${\cal D}_{*}$ satisfies condition  \e{eq:D}. Let the measures $u(z) dM(z)$ be absolutely continuous
      for all $u\in{\cal D}_{*}$.  Then the measure $  dM(z)$  is also  absolutely continuous.
    \end{lemma}

 Now we are in a  position to conclude the proof of  Theorem~\ref{T1}. 
  Suppose that the form $t[f,f] $  is closable. Then by Lemma~\ref{adj} the condition  \e{eq:D} is satisfied. 
  By the definition of the set ${\cal D}_{*}$, the Fourier coefficients of  the measures $d \mu(z)= u(z) dM(z)$ belong to $\ell^2 ({\Bbb Z}_{+})$ for all $u\in{\cal D}_{*}$. Therefore it follows from   Theorem~\ref{brothers} that   these measures   are absolutely continuous.  
        Hence   by Lemma~\ref{ac}, the measure $  dM(z)$  is also absolutely continuous.



  \medskip
  
  {\bf 2.5.}
Next, we sketch  the proof of Theorem~\ref{Hamb1}. It is almost obvious that the conditions $(ii)$ and $ (iii)$ are equivalent. So we discuss only the equivalence of $(i)$ and $ (ii)$. Algebraically, we follow the scheme of the previous section, but instead of $L^2 (\Bbb T; dM)$ we introduce the space $L^2 (\Bbb R; d{\sf M})$ where the measure $d{\sf M}$ is defined on $\Bbb R$ by equations \e{eq:WHh}.  
The role  of the operator $A$ is now played by the operator $B \colon \ell^2 ({\Bbb Z}_{+})\to   L^2 (\Bbb R; d{\sf M})$  defined   on the set $\cal D$ by the formula
   \begin{equation}
( Bf) (x)=\sum_{n=0}^\infty   f_{n} x^n , \q x\in {\Bbb R}.
 \label{eq:AH}\end{equation}
  Instead of \e{eq:A2}, we now have the identity
  \begin{equation}
    \| B f\|^2_{L^2 (\Bbb R; d{\sf M})} = h[f,f] , \q f\in \cal D,
 \label{eq:AHx}\end{equation}
  and the role of Lemma~\ref{de} is played by the following assertion.

    \begin{lemma}\label{deH}
 The form $   h[f,f] $ defined   on $\cal D$   is closable  in the space $\ell^2 ({\Bbb Z}_{+})$  if and only if the operator $B $   is closable.  
    \end{lemma}
    
    The adjoint operator  
    $B^* \colon L^2 (\Bbb R; d{\sf M}) \to \ell^2 ({\Bbb Z}_{+}) $
    can be constructed similarly to Lemma~\ref{LTM}.

   \begin{lemma}\label{LTMH}
   Let a subset ${\sf D}_{*}$ of $L^2 (\Bbb R; d{\sf M})$ consist of functions $u(x)$ such that the sequence
    \[
u_{n}:= \int_{-\infty}^\infty  u(x) x^n d{\sf M}(x) 
 \]
 belongs to $\ell^2 ({\Bbb Z}_{+})$. Then  the operator $B^*$ is given by the equality 
$      (B^ {*}u)_{n}=u_{n} $
  on the domain ${\cal D}(B^*) ={\sf D}_{*}$. 
    \end{lemma}
    
 For detailed proofs of these assertions see Lemmas~2.1 and 2.2 in \cite{Yunb}. 

 The ``if " part of   the following  result is quite easy, but the converse statement requires the Paley-Wiener theorem.

     
     \begin{theorem}\label{dense}\cite[Theorem~2.3]{Yunb} 
     The set ${\sf D}_{*}$ is dense in $L^2 ({\Bbb R}; d{\sf M})$ if and only if condition \e{eq:dense} is satisfied.
    \end{theorem}
    
   We  only make some comments on the proof  of the `` only if " part. Actually, only the inclusion
        \begin{equation}
 \supp {\sf M} \subset [-1,1]
\label{eq:X7}\end{equation}
deserves a special discussion.
    
    
     For an arbitrary $u\in L^2({\Bbb R}; d{\sf M})$, we put
 \begin{equation}
\Psi (z) = \int_{-\infty}^\infty e^{i z x} u(x) d{\sf M} (x).
\label{eq:X1}\end{equation}
Since all functions $x^n$ belong to $ L^2({\Bbb R}; d{\sf M})$, we see that
$\Psi \in C^\infty ({\Bbb R})$ and
 \[
\Psi^{(n)} (0)= i^n \int_{-\infty}^\infty  x^n u(x) d{\sf M}(x).
\]
If $u\in{\sf D}_{*}$,
then this sequence is bounded and hence the function
\[
\Psi(z)= \sum_{n=0}^\infty \frac{\Psi^{(n)} (0)}{n!} z^n 
\]
is entire and satisfies the estimate
\[
|\Psi (z)|\leq C\sum_{n=0}^\infty \frac{1}{n!} |z|^n =C e^{|z|}, \q z\in {\Bbb C}.
\]
By virtue of the Phragm\'en-Lindel\"of principle, we actually have a stronger estimate
\begin{equation}
|\Psi(z)|\leq  C e^{ |\Im z|}, \q z\in {\Bbb C}.
\label{eq:X5}\end{equation}

According to the Paley-Wiener theorem (see, e.g., Theorem~IX.12 in \cite{RS}) it follows from estimate \e{eq:X5} that the Fourier transform of $\Psi (z)$ (considered as a distribution in the Schwartz class ${\cal S}' (\Bbb R)$) is supported by the interval $[ -1,1]$. Therefore  formula \e{eq:X1} implies that  
 \begin{equation}
 \int_{-\infty}^\infty \varphi(x)u(x) d {\sf M}(x) =0, \q \forall \varphi\in C_{0}^\infty ({\Bbb R}\setminus [-1,1]),
\label{eq:X6}\end{equation}
for all $ u\in{\sf D}_{*}$.
If ${\sf D}_{*}$ is   dense in $L^2 ({\Bbb R}; d{\sf M})$, then we can approximate $1$ by functions  $u\in {\sf D}_{*}$  in this space. 
Hence equality \e{eq:X6} is true with $u(x)=1$ which implies \e{eq:X7}.    $\Box$

Finally, we almost repeat the arguments used in the proof of Theorem~\ref{T1}.  Putting together Lemma~\ref{LTMH} and Theorem~\ref{dense}, we see that the operator $B^*$ is densely defined and hence $B$  is closable if and only if condition \e{eq:dense} is satisfied. In view of Lemma~\ref{deH} this proves that the conditions (i) and (ii) of Theorem~\ref{Hamb1} are equivalent.  $\Box$ 


    \medskip
       

{\bf 2.6.}
According to Theorem~\ref{Hamb1}, the condition        $ h_n\to 0$ as $n\to \infty$ is necessary  and sufficient for a   Hankel quadratic form \e{eq:QFh} to be   closable.   On the contrary, it is probably impossible to give necessary and sufficient conditions (at least elementary) for a Toeplitz quadratic form \e{eq:QFq} to be closable  in terms of its  entries $t_{n}$. Indeed, an obvious necessary condition is $t_{n}\to 0$  as $| n| \to\infty$ because, by Theorem~\ref{T1},  the measure $dM(z)$ in the representation   \e{eq:WH}  is absolutely continuous. An obvious sufficient condition is $ \{t_n\} \in \ell^2 ({\Bbb Z})$ because in this case, by the Parseval identity,  the measure $dM(z)$ is absolutely continuous and its derivative $t\in L^2 ({\Bbb T} )$.

            Apparently,  this      gap between necessary  and sufficient conditions  cannot be significantly reduced. Indeed, by the Wiener theorem (see, e.g., Theorem~XI.114 in \cite{RS}), if the Fourier coefficients $ t_n$ of some measure $dM(z)$  tend to zero, then this measure is necessarily continuous, but it may be singular with respect to the Lebesgue measure.   Thus  the condition        $ t_n\to 0$ as $|n |\to \infty$ does not imply that the   measure $dM(z)$ defined by equations \e{eq:WH} is absolutely continuous. So in accordance with  Theorem~\ref{T1}    
   the   corresponding Toeplitz quadratic form $t[f,f]$ need not be closable.
        
          Astonishingly, the   condition  $ \{t_n\}_{n\in {\Bbb Z}}\in \ell^2 ({\Bbb Z})$  guaranteeing   the absolute continuity of the measure $dM(z)$ turns out to be very sharp.  Indeed,  for every $p \in{\Bbb Z}_{+}$, O.~S.~Iva\v{s}\"{e}v-Musatov constructed in \cite{I-M}  a singular measure such that its Fourier coefficients satisfy the estimate
  \[
  t_{n} =O \big( ( n (\ln n)( \ln\ln n)\cdots (\ln_{(p)} n))^{-1/2}  \big)
  \]
  (here $ \ln_{(p)} n$ means that the logarithm  is applied $p$ times to $n$). This sequence ``almost belongs" to $\ell^2 ({\Bbb Z})$, but, by Theorem~\ref{T1}, the   corresponding   form $t[f,f]$ is not  closable.


   
   \section{Closed quadratic forms} 

 Here we will show that  closable Toeplitz $t[f,f]$ and Hankel $h[f,f]$  quadratic forms constructed in Theorems~\ref{T1} and \ref{Hamb1}, respectively, are closed on their maximal domains of definition.
This yields a  description  of the domains of the operators ${\cal D} (\sqrt{T})$ and  ${\cal D} (\sqrt{H})$.  Again, the algebraic scheme of a study of  Toeplitz   and Hankel   forms is the same, but analytical details are quite different. As before, it is convenient to use the operators $\bf A$ and $B$ defined by formulas \e{eq:A} and \e{eq:AH}.

  \medskip
  
{\bf 3.1.}
Let us start with  closable Toeplitz forms when, by Theorem~\ref{T1}, the measures $dM(z)$ in equations \e{eq:WH} are absolutely continuous, that is,   $dM(z)= t(z) d{\bf m} (z)$ where $t\in L^1 ({\Bbb T} )$ and
we may suppose that    $t(z)\geq 1$. 

 
 Under the assumptions of Theorem~\ref{T1} the  operator $A^*$ adjoint to $A$ (recall that $A$ is the restriction of $\bf A$ of $\cal D$)  is densely defined so that the second adjoint exists and $A^{**}=\cl A$ (the closure of $A$).   Let us also  introduce by formula \e{eq:A} the ``maximal" operator $A_{\rm max}$   on the domain 
${\cal D} (A_{\rm max})$ that consists of all $f\in \ell^2 ({\Bbb Z}_{+})$ such that ${\bf A}f\in L^2 (\Bbb T; dM)$.
We will show that
          \begin{equation}
     \cl{A}= A_{\rm max}.
     \label{eq:MAX}\end{equation}
    
    The first assertion is a direct consequence of the definition of the closure   of the operator $A$.
    
     \begin{lemma}\label{Yatr} 
     The inclusion $\cl{A}\subset A_{\rm max}$ holds.
        \end{lemma}
        
        Indeed, if $f\in {\cal D}(\cl{A})$, 
        then there exists a sequence $f^{(k)}\in {\cal D}( A)$ such that $f^{(k)}\to f$
         in $\ell^2 ({\Bbb Z}_{+})$ and $ A f^{(k)} $ converges to some $ g\in  L^2 (\Bbb T; dM)$; 
         in this case, $(\cl A ) f:=g$. It follows from \e{eq:A2P} that
        ${\bf A}f^{(k)}\to {\bf A}f$ in $L^2 (\Bbb T )$ whence $g= {\bf A}f$.  Since $ g\in L^2 (\Bbb T; dM)$,  we have $g=A_{\rm max}f$ where $f\in {\cal D}(A_{\rm max})$.

    The opposite inclusion   is less trivial.
 
 \begin{lemma}\label{YaTL} 
 The inclusion $   A_{\rm max}\subset \cl{A}$ holds.
        \end{lemma}
        
         Pick  $ f\in{\cal D} ( A_{\rm max})$. Then $f\in \ell^2 ({\Bbb Z}_{+})$ and   $u: ={\bf A} f\in H^2 (\Bbb T)\cap L^2 (\Bbb T; dM)  $. Since $t \in L^1(\Bbb T)$ and $t(\z
         )\geq 1$,  
         \[
   t_{\rm out}(z) :=      \exp\Big( \frac{1}{2}\int_{\Bbb T} \frac{\z +z}{ \z -z}\ln t(\z) d{\bf m} (\z) \Big),\q |z|<1,
         \]
            is an           outer function,        $t_{\rm out}\in H^2 (\Bbb T)$ and the angular limits of $t_{\rm out}(z)$ as $z\to \z\in{\Bbb T}$   equal        $\sqrt{t(\z)}$.  We also have $v: =u t_{\rm out} \in H^2 (\Bbb T)$ because   $u\sqrt{t} \in               
            L^2 (\Bbb T)  $.
                  By the V.~Smirnov theorem (see, e.g., \cite{NK}, Section1.7), every function in $ H^2 (\Bbb T)$ can be approximated by linear combinations of functions $z^n t_{\rm out}(z)$ so that there exists a sequence of polynomials ${\cal P}^{(k)} (z)$ such that
                  \[
            \| v- {\cal P}^{(k)}t_{\rm out}\|_{ H^2 (\Bbb T)}^2 =\int_{\Bbb T}| u(z) - {\cal P}^{(k)} (z)|^2 t(z)d{\bf m}  (z)\to 0
          \]
      as $   k\to\infty  $.    This means that
          \begin{equation}
\lim_{k\to\infty} \| {\bf A}f  - {\bf A}f^{(k)} \|_{L^2 (\Bbb T; dM)}=0
 \label{eq:BC}\end{equation}
 for  $f^{(k)}\in {\cal D}$ such that ${\bf A}f^{(k)}= {\cal P}^{(k)}$.
  Since the convergence in $ L^2 (\Bbb T; d M)$ is stronger than   in $ L^2 (\Bbb T )$, we see that ${\bf A}f^{(k)} \to 
{\bf A}f$  in $ L^2 (\Bbb T )$ as $k\to\infty$. According to    \e{eq:A2P}  this     implies that 
 $ f^{(k)} \to 
 f$  in $  \ell^2 ({\Bbb Z}_{+})$ as $k\to\infty$. Combining this relation with \e{eq:BC}, we see that $ f\in{\cal D} ( \cl A )$ and $u=( \cl A)f. \q \Box$ 
         
        In view of relations   \e{eq:A2P} and \e{eq:A2}, equality  \e{eq:MAX} can be reformulated in terms of Toeplitz quadratic forms $t[f,f]$.


\begin{theorem}\label{Clo}\cite[Theorem~2.10]{YaT} 
 Under the assumptions of Theorem~\ref{T1} 
the closure  of the form $t[f,f]$ is given by the equality 
  \begin{equation}
   t [f,f]     =\int_{\Bbb T}  | ({\bf A} f) (z)|^2 dM(z)  
 \label{eq:A2x}\end{equation}
  on the set ${\cal D}[t] = {\cal D} (A_{\rm max})$ of all $f\in \ell^2 ({\Bbb Z}_{+})$ such that the right-hand side of \e{eq:A2x} is finite.
        \end{theorem}


 \medskip  
  
{\bf 3.2.}
For Hankel quadratic forms, we proceed  from Theorem~\ref{Hamb1}. We suppose that  condition \e{eq:dense}  is satisfied, and hence  the form $h[f,f]$ is closable. Let us  now define the operator ${\bf B} $ by formula \e{eq:AH} for all $f\in \ell^2({\Bbb Z}_{+})$. The series  in the right-hand side of  \e{eq:AH} converges for   each  $x\in(-1,1)$, but only the estimate
 \[
|({\bf B} f)(x)|\leq \sum_{n=0}^\infty |f_{n}| |x|^n \leq  (1-x^2)^{-1/2}\,  \| f\|_{\ell^2 ({\Bbb Z}_{+})} 
 \]
 holds. So it is of course possible that ${\bf B}  f \not\in L^2 ((-1,1); d{\sf M})$. Therefore we also introduce the ``maximal" operator $B_{\rm max}$ as the restriction of ${\bf B}$ on the domain 
${\cal D} (B_{\rm max})$ that consists of all $f\in \ell^2 ({\Bbb Z}_{+})$ such that ${\bf B} f\in L^2 ((-1,1); d{\sf M})$. The following result plays the role of equality \e{eq:MAX}.

 \begin{lemma}\label{YaHL} 
     Let  one of equivalent conditions of Theorem~\ref{Hamb1} be satisfied.    
       Then  equality 
        \begin{equation}
     \cl{B}= B_{\rm max}.
    \label{eq:MAH}\end{equation} 
     is true.
        \end{lemma}
         
             Similarly to Lemma~\ref{Yatr}, the inclusion $\cl{B}\subset B_{\rm max}$ is a direct consequence of the definition of the closure  of the operator $B$. Surprisingly, the proof of the opposite inclusion 
             $B_{\rm max}  \subset B^{**} $  turns out to be rather tricky although it does not require any deep analytical results. Of course, it suffices to check that
          \[
     (B_{\rm max}f,u)_{L^2 ((-1,1); d{\sf M})} = (f, B^{*} u)_{\ell^2 ({\Bbb Z}_{+})}
 \]
 for all $f\in {\cal D} (B_{\rm max})$ and all $u\in {\cal D} (B^*)= {\sf D}_{*}$. In the detailed notation, this relation   means that
           \begin{equation}
  \int_{-1}^1 \big(\sum_{n=0}^\infty f_{n}  x^n  \big) \ov{u(x)} d{\sf M}(x)
  = \sum_{n=0}^\infty f_{n}\big(   \int_{-1}^1   x^n    \ov{u(x)} d{\sf M}(x) \big). 
\label{eq:AB4}\end{equation}
 The problem is that these integrals do not converge absolutely. So the Fubini theorem cannot be applied, and we have not found a direct proof of  relation \e{eq:AB4}.
 
  By some, rather mysterious reasons, it appears to be more convenient to treat this problem in the realization of Hankel operators as integral operators in the space $L^2({\Bbb R}_{+})$. This means that instead of the operator $\bf B$ defined by formula \e{eq:AH} for all $f\in \ell^2 ({\Bbb Z}_{+})$, we now consider the operator (the Laplace transform)
 defined by the formula
   \begin{equation}
  ({\bf G}f) (\lambda)= \int_{0}^\infty e^{-t \lambda} f(t) dt  
\label{eq:LAPj}\end{equation}
for all $f\in L^2 ({\Bbb R}_{+})$. The role of $L^2 ((-1,1); d{\sf M})$ is played by the space $L^2 ({\Bbb R}_{+}; d\Sigma)  $ 
 where the non-negative measure $ d\Sigma (\lambda)$ on ${\Bbb R}_{+}$ satisfies the condition 
  \begin{equation}
\int_{0}^\infty (\lambda +1 )^{-k} d\Sigma (\lambda) <\infty 
 \label{eq:SS}\end{equation}
  for   $k=2$.         The integral \e{eq:LAPj} converges for all $f\in L^2 ({\Bbb R}_{+})$ and $\lambda>0$, but  the estimate 
 \[
 |   ({\bf G}  f) (\lambda)|\leq (2\lambda)^{-1/2} \, \| f\|_{L^2 ({\Bbb R}_{+})}  
 \]
 does not of course  guarantee that $ {\bf G} f\in L^2 ({\Bbb R}_{+}; d\Sigma)$.


  Properties of the operators ${\bf B}$ and ${\bf G}$ are basically the same.
   We first define the restriction $G$ of the operator  ${\bf G}$  on domain   ${\cal D}(G)$
that consists of functions compactly supported in ${\Bbb R}_{+}$. Evidently, $Gf\in  L^2 ({\Bbb R}_{+}; d\Sigma)$ if
$f\in {\cal D}(G)$. 
    It is easy to show (see \cite{Yf1a}, for details) that the operator $G^*$ is given by the formula
 \[
  (G^* v) (t)= \int_{0}^\infty e^{-t \lambda} v(\lambda) d\Sigma(\lambda) ,
\]
and $v\in {\cal D} (G^*)$ if and only if $v\in L^2 ({\Bbb R}_{+}; d\Sigma)$ and $G^* v\in L^2 ({\Bbb R}_{+})$.
Obviously, this condition is satisfied if $v$ is compactly supported in ${\Bbb R}_{+}$.
Since the set of such $v$  is dense in $L^2 ({\Bbb R}_{+}; d\Sigma)$, the operator $G^*$ is densely defined. Thus $G$ is closable   and $\cl G=G^{**}$.
Let us now define
   the operator $G_{\rm max}$   as the restriction of the operator ${\bf G}$  on the domain $\cal D(G_{\rm max})$ that consists of all $f\in L^2 ({\Bbb R}_{+})$ such that ${\bf G} f\in L^2 ({\Bbb R}_{+}; d\Sigma)$.

  The following assertion plays the central role.

\begin{lemma}\cite[Theorem~3.9]{Yf1a}\label{LLL}  
 Let     $ d\Sigma (\lambda)$ be a measure  on ${\Bbb R}_{+}$ such that the condition 
 \e{eq:SS} is satisfied for some $k>0$. Then
       \begin{equation}
\cl G  = G_{\rm max}.
 \label{eq:ABL}\end{equation}
        \end{lemma}
        
        We will not comment on a rather complicated proof of this result, but explain the equivalence of relations \e{eq:MAH} and 
        \e{eq:ABL} (for the particular case $k=2$). 
Suppose that the   measures $ d\Sigma (\lambda)$ and $d {\sf M}(x)$ are linked by the equality
  \[
d {\sf M}(x)=(\lambda +1/2)^{-2} d\Sigma (\lambda), \q x=\frac{2\lambda-1}{2\lambda+1}.
 \]
 Thus,  ${\sf M}((-1,1))<\infty$ if and only if the condition 
 \e{eq:SS} holds for $k=2$.
Let us also set 
 \[
  (Vu) (\lambda)= \frac{1}{\lambda+1/2} u \Big( \frac{2 \lambda-1}{2\lambda+1} \Big).  
\]
Obviously, $V: L^2 ((-1,1); d{\sf M}) \to L^2 ({\Bbb R}_{+}; d\Sigma)$ is a unitary operator.

We need  the identity (see formula (10.12.32) in \cite{BE})
 \begin{equation}
  \int_{0}^\infty     {\sf L}_{n} (t ) e^{-(1/2+ \lambda)t} dt=  \frac{1}{\lambda+1/2}\Big( \frac{2 \lambda-1}{2\lambda+1} \Big)^n, \q \lambda> -1/2 , 
\label{eq:KL4x}\end{equation}
for the Laguerre polynomials ${\sf L}_{n} (t )$ (see, for example,  the book \cite{BE}, Chapter~10.12, for their definition). 
It can be deduced from this fact that the functions 
$  {\sf L}_{n}  (t) e^{-t/2}$, $n=0,1,\ldots$,   
form an orthonormal basis in the space $ L^2({\Bbb R}_{+})$, and hence the operator  $U \colon l^2({\Bbb Z}_{+})\to  L^2({\Bbb R}_{+})$ defined by the formula  
\begin{equation}
(U  f ) (t)=\sum_{n=0}^\infty f_{n}  {\sf L}_{n}  (t) e^{-t/2}, \q f=(f_{0}, f_{1}, \ldots),
\label{eq:K3}\end{equation}
is unitary.

A link of the operators $B$ and $G$ is stated in the following assertion which can be easily derived from \e{eq:LAPj}, \e{eq:KL4x}   and \e{eq:K3}; see Lemma~3.2 in \cite{Yunb}, for details. 

 \begin{lemma}\label{LA} 
     For all $f\in{\cal D}$, the identity  holds
      \[
V {\bf B}  f= {\bf G} U  f.
\]
        \end{lemma}


 Combining    Lemmas~\ref{LLL} and \ref{LA}, we arrive at
 the following result.
 
 \begin{lemma}\label{LM}  
 Let     $ d {\sf M} (x)$ be a finite measure  on $(-1,1)$. Then    equality \e{eq:MAH} holds.
        \end{lemma}
        


  In view of  identity \e{eq:AHx}, equality \e{eq:MAH}  leads to the following result which plays the role of Theorem~\ref{Clo}.

  \begin{theorem}\label{ClH}\cite[Theorem~3.4]{Yunb} 
     Let  the form $h[f,f]$ be defined on  ${\cal D}$ by  \e{eq:Ha}, and let   assumption \e{eq:QFh} be true. Suppose that one of three equivalent conditions (i), (ii) or (iii) of Theorem~\ref{Hamb1} is satisfied. 
   Then  the closure of $h[f,f]$
 is given by the equality 
 \begin{equation}
   h[f,f] =  \int_{-1}^1  | \sum_{n=0}^\infty   f_{n} x^n|^2 d {\sf M}(x) 
 \label{eq:AB2}\end{equation}
  on the set  of all $f\in \ell^2 ({\Bbb Z}_{+})$ such that the right-hand side of \e{eq:AB2} is finite.
        \end{theorem}

Note that that domains ${\cal D} (T)$ and  ${\cal D} (H)$ of Toeplitz $T$ and Hankel $H$ do not admit an explicit description, but
 Theorems~\ref{Clo} and \ref{ClH}  characterize the domains 
  ${\cal D} (\sqrt{T})$ and  ${\cal D} (\sqrt{H})$ of their square roots.
  
     \section{Wiener-Hopf semibounded operators} 
     
 {\bf 4.1.}
Wiener-Hopf operators   $W$ are formally  defined in the space $L^2({\Bbb R}_{+})$ of functions  $f(x)$
 by the formula
  \[
(W f) (x)=\int_{ {\Bbb R}_{+}}   w(x-y) f(y) dy.
\]
 These operators are continuous analogues of the Toeplitz operators defined by \e{eq:Toe}.
  However optimal results on 
 Wiener-Hopf  operators are not direct consequences of the corresponding results for Toeplitz operators and, in some sense, they are more general.  One of the differences is that  Wiener-Hopf  operators  require a consistent work with distributions.
 
To be precise, we define the operator $W$ via its quadratic form
   \begin{equation}
w[f,f] =\int_{{\Bbb R}_{+}}  \int_{{\Bbb R}_{+}}     w(x-y) f(y) \ov{f (x)}dx dy .
 \label{eq:WFq}\end{equation}
  With respect to   $w$, we a priori only assume that it is a distribution in the class $  C_{0}^\infty ({\Bbb R})'$ dual to $  C_{0}^\infty ({\Bbb R})$. Then the quadratic form is correctly defined for all  $f\in C_{0}^\infty ({\Bbb R}_{+})$.
   We always suppose that $w(x)= \ov{w(-x)}$ so that the operator $W$ is formally symmetric and the quadratic form \e{eq:WFq} is real. We also assume that it is semibounded from below, that is,
    \begin{equation}
w[f,f] \geq \gamma \| f\|^2  ,\q f\in C_{0}^\infty ({\Bbb R}_{+}), \q \| f\| =\| f\|_{L^2({\Bbb R}_{+})},
 \label{eq:Tw1}\end{equation}
 for some  $\gamma\in \Bbb R$.

Here we follow basically the   scheme  we used before for semibounded  Toeplitz operators. However the analytical basis is   rather different.  
The role of Theorem~\ref{R-H} is now played by  the Bochner-Schwartz theorem (see, e.g., Theorem~3 in \S 3 of Chapter II of the book \cite{GUEVI+}).

  \begin{theorem}\label{Bochner} 
  Let the form $w[f,f] $ be defined by the relation \e{eq:QFq} where the distribution $w\in C_{0}^\infty({\Bbb R})'$.  Then the condition
      \[
w [f,f]  
 \geq 0, \q \forall f\in C_{0}^\infty ({\Bbb R}_{+}),
\]
  is satisfied if and only if there exists a
  non-negative   measure $d {\rm M} (\lambda)$ on the line $\Bbb R$ such that 
  \begin{equation}
w(x) = \frac{1}{2\pi} \int_{\Bbb R}  e^{-ix \lambda} d{\rm M}(\lambda ).
 \label{eq:WHW}\end{equation}
Here the measure obeys the condition
    \begin{equation}
\int_{\Bbb R} (1+\lambda^2)^{-p} d{\rm M}( \lambda)<\infty
\label{eq:Sch}\end{equation}
for some $p $ $($that is,  it has at most a polynomial growth at infinity$)$.
    \end{theorem}

Observe that the Lebesgue measure $d{\rm M}( \lambda)=d\lambda$ satisfies the condition \e{eq:Sch} with $p>1/2$. 
   For  the Lebesgue measure,  relation \e{eq:WHW} yields $w(x)=\d(x)$ (the delta-function) so that  $W=I$ and  $w[f,f] =  \| f\|^2$.  Therefore the measure corresponding to the form
    $w[f,f]+ \beta \| f\|^2 $ equals $d{\rm M}(\lambda ) + \beta d \lambda$, and  relation \e{eq:WHW} extends to all semibounded Wiener-Hopf quadratic forms.   Thus we have the one-to-one correspondence between Wiener-Hopf  quadratic forms satisfying estimate  \e{eq:Tw1} and real measures satisfying the condition $ {\rm M}  (X)\geq \gamma    |X|$ ($  |X|$ is the Lebesgue measure of $X $) for all Borelian sets $X\subset{\Bbb R}$.
   
   \medskip

{\bf 4.2.} 
  Our goal is to find necessary and sufficient conditions for the form $w[f,f]$ to be  closable. The following result plays the role   of Theorem~\ref{T1}.
   
    \begin{theorem}\label{TW1} 
    Let the form $w[f,f]$ be given by formula \e{eq:WFq} on elements $f\in C_{0}^\infty ({\Bbb R}_{+})$, and let the condition 
 \e{eq:Tw1} be satisfied  for some  $\gamma\in \Bbb R$.
   Then the form $w[f,f] $ is closable in the space $L^2 ({\Bbb R}_{+})$ if and only if  the measure $d{\rm M}  (\lambda)$ in the equation \e{eq:WHW} is  absolutely continuous.
    \end{theorem}
    
 We always understand the absolute continuity with respect to the Lebesgue measure.   Therefore  Theorem~\ref{TW1} means that $d{\rm M}  (\lambda)=  \varphi (\lambda) d\lambda$ where $\varphi \in L^1_{\rm loc} (\Bbb R)$,
   \[
\int_{\Bbb R} (1+\lambda^2)^{-p} | \varphi (\lambda)|d  \lambda<\infty
 \]
       and $\varphi (\lambda) \geq \gamma$.  
       The function $\varphi(\lambda)$ is known as the symbol of the Wiener-Hopf operator $W$. Thus Theorem~\ref{TW1} shows that in the semibounded case, the symbol of a  Wiener-Hopf operator can be correctly defined if and only if the corresponding quadratic form is closable.

 Our proof of Theorem~\ref{TW1}   requires a
continuous analogue of the classical Riesz Brothers theorem. Let us state this result here. For a measure $d {\rm M} ( \lambda)$ on $\Bbb R$, we denote by $d | {\rm M} | ( \lambda)$ its  variation.

     \begin{theorem}\label{RBT} 
     Let   $d {\rm M} (\lambda)$ be a complex     measure on the line $\Bbb R$ such that
      \[
\int_{\Bbb R} (1+\lambda^2)^{-p} d | {\rm M} | ( \lambda)<\infty
\]
for some $p$. Put
  \[
\sigma(x) = \frac{1}{2\pi} \int_{\Bbb R} e^{-ix \lambda} d{\rm M} (\lambda ) 
\]
 and suppose that $\sigma \in L^2 (a,\infty)$ for some $a\in \Bbb R$. Then the measure $d {\rm M} (\lambda)$  is  absolutely continuous.
    \end{theorem}   
    
    We allow $a\in \Bbb R$ in   Theorem~\ref{RBT} to be arbitrary since, for example, the function $\sigma(x)=\d (x-x_{0})$ for any $x_{0}\in \Bbb R$ does not belong to $L^2_{\rm loc} ({\Bbb R})$, but the corresponding measure $d {\rm M} (\lambda)= e^{i x_{0}\lambda} d\lambda$ is of course absolutely continuous.

 \medskip

  {\it Acknowledgements}. The author thanks N.~K.~Nikolski for a useful discussion. 



\begin{thebibliography}{00}

\bibitem{AKH} N.~Akhiezer,
 \emph{The classical moment problem and some related questions in analysis},  Oliver and  Boyd, Edinburgh and London, 1965.

 \bibitem{BSbook}
  M.~Sh.~Birman, M.~Z.~Solomyak, 
 \emph{Spectral theory of selfadjoint operators in Hilbert space}, D.~Reidel, Dordrecht, 1987. 

\bibitem{Bo} A.~B\"ottcher, B.~Silbermann,
\emph{Analysis of Toeplitz operators},  Springer-Verlag,  2006.

 
\bibitem{Carleman} T.~Carleman,
\emph{Sur les \'equations int\'egrales singuli\`eres},  Uppsala,  1923.



\bibitem{BE} A.~Erd\'elyi,  W.~Magnus, F.~Oberhettinger, F.~G.~Tricomi,
\emph{Higher transcendental functions}, Vol. 1, 2, 
McGraw-Hill, New York-Toronto-London, 1953.



\bibitem{GUEVI+} I.~M.~Gel'fand,  G.~E.~Shilov,  
\emph{Generalized functions.  } Vol.~1, 
Academic Press, New York and London,  1964.   

\bibitem{GGK}
 I.~Gohberg, S.~Goldberg, M.~Kaashoek,
\emph{Classes of linear operators},   vol. 1, Birkh\"auser,  1990. 

 \bibitem{Hamb}
 H.~Hamburger,
\emph{\"Uber eine Erweiterung des Stieltjesschen Momentenproblems}, Math. Ann.  {\bf 81} (1920), 235-319;  {\bf 82}  (1921), 120-164 and 168-187.


\bibitem{Hartman}
P.~Hartman,  
\emph{On unbounded Toeplitz matrices}, Amer. J. Math.   {\bf 85}, N 1  (1963), 59-78.   



\bibitem {Hof} K. Hoffman, {\it Banach spaces of analytic functions}, Prentice-Hall, Inc., Englewood Cliffs, New York, 1962.

  \bibitem{I-M}
O.~S.~Iva\v{s}\"{e}v-Musatov,  
\emph{On the Fourier-Stieltjes coefficients of singular functions},  Dokl. Akad. Nauk SSSR   {\bf 82}  (1952), 9-11 (Russian).


 


 
\bibitem{Nehari}
Z.~Nehari,  
\emph{On bounded bilinear forms}, Ann. Math.   {\bf 65} (1957), 153-162. 
 

 


 \bibitem{NK} N.~K.~Nikolski, {\em Operators, functions, and systems: an easy reading}, vol. I: Hardy, Hankel, and Toeplitz, Math. Surveys and Monographs vol.~92,  Amer. Math. Soc.,   Providence,
  Rhode Island, 2002.
  
 
  
\bibitem{Pe}
V.~V.~Peller,
\emph{Hankel operators and their applications}, 
Springer Verlag,  2002. 

\bibitem {RS} M. Reed, B. Simon, {\em Methods of modern mathematical physics}, vol. 2,
 Academic Press, San Diego, CA,  1975.

   

\bibitem{Sarason}
D.~Sarason, \emph{Unbounded Toeplitz operators},  Integral Eq. Op. Th.   {\bf 61} (2008), 281-298.



\bibitem{Widom}
H.~Widom,
\emph{Hankel matrices}, 
Trans. Amer. Math. Soc. {\bf 121}  (1966), 1-35.


  
 


 



\bibitem{Yf1a}
D.~R.~Yafaev, \emph{Quasi-Carleman operators and their spectral properties},  Integral Eq. Op. Th.   {\bf 81} (2015), 499-534.

\bibitem{Yunb}
D.~R.~Yafaev, \emph{ Unbounded Hankel operator and moment problems},   Integral Eq. Oper. Theory, {\bf 85}  (2016), 289-300.

\bibitem{YaT}
 D.~R.~Yafaev, \emph{ On semibounded Toeplitz operators},   J. Operator Theory,  {\bf 77}  (2017), 742-762.

\bibitem{YaWH}
D.~R.~Yafaev, \emph{ On semibounded Wiener-Hopf operators}, Journal of LMS,  {\bf 95}  (2017), 101-112.



\end{thebibliography}
 \end{document}